\documentclass[leqno,oneside,english]{article}
\usepackage[T1]{fontenc}
\usepackage[latin1]{inputenc}
\usepackage{amsmath}
\usepackage{amssymb}

\newcommand{\findem}{{\unskip \nobreak \hfil \penalty 50\hskip 2em\hbox {}\nobreak
\hfil \pave \parfillskip =0pt\finalhyphendemerits =0\smallbreak }}
\newcommand{\pave}{\hbox {\vrule height5pt width5pt}}
\def\Vo{\vbox{\offinterlineskip\hbox{\kern 3pt$\scriptstyle\circ$}
\kern 1pt\hbox{$V$}}}

\def\Ho{\vbox{\offinterlineskip\hbox{\kern 3pt$\scriptstyle\circ$}
\kern 1pt\hbox{$H$}}}

\def\Wo{\vbox{\offinterlineskip\hbox{\kern 3pt$\scriptstyle\circ$}
\kern 1pt\hbox{$W$}}}

\def\and{\quad\mbox{and}\quad}

\def\rot{\mathop{\bf curl}\nolimits}

\newcommand{\be}{\begin{equation}}
\newcommand{\ee}{\end{equation}}
\newcommand{\beq}{\begin{eqnarray}}
\newcommand{\eeq}{\end{eqnarray}}
\newcommand{\beqs}{\begin{eqnarray*}}
\newcommand{\eeqs}{\end{eqnarray*}}
\newcommand{\bt}{\begin{theorem}}
\newcommand{\et}{\end{theorem}}
\newcommand{\br}{\begin{remark}}
\newcommand{\er}{\end{remark}}
\newcommand{\bc}{\begin{corollary}}
\newcommand{\ec}{\end{corollary}}
\newcommand{\bl}{\begin{lemma}}
\newcommand{\el}{\end{lemma}}
\newcommand{\bp}{\begin{proposition}}
\newcommand{\ep}{\end{proposition}}
\newcommand{\bd}{\begin{definition}}
\newcommand{\ed}{\end{definition}}
\newcommand{\bex}{\begin{example}}
\newcommand{\eex}{\end{example}}
\newcommand{\pf}{{\bf Proof: }}

\newcommand{\C}{{\mathbb C}}
\newcommand{\R}{{\mathbb R}}

\newcommand{\eps}{\varepsilon }

\newcommand{\caD}{{\cal D}}

\newcommand{\caH}{{\cal H}}

\newcommand{\ddiv}{{\rm div}}
\newcommand{\Om}{{\Omega}}
\newtheorem{theorem}{Theorem}[section]
\newtheorem{lemma}[theorem]{Lemma}
\newtheorem{corollary}[theorem]{Corollary}
\newtheorem{remark}[theorem]{Remark}
\newtheorem{definition}[theorem]{Definition}
\newtheorem{proposition}[theorem]{Proposition}
\newtheorem{example}[theorem]{Example}

\def\be{\begin{equation}}
\def\ee{\end{equation}}
\def\bea{\begin{eqnarray}}
\def\eea{\end{eqnarray}}
\def\ea{\end{array}}
\def\ds{\displaystyle}

\def\rot{{\rm curl}}

\def\11{{\rm 1~\hspace{-1.5ex}1} }

\def\CC{\rm \hbox{C\kern-.56em\raise.4ex
         \hbox{$\scriptscriptstyle |$}\kern+0.5 em }}

\newcommand{\rfb}[1]{\mbox{\rm
   (\ref{#1})}\ifx\undefined\stillediting\else:\fbox{$#1$}\fi}

\makeatletter
\def\section{\@startsection {section}{1}{\z@}{-3.5ex plus -1ex minus
    -.2ex}{2.3ex plus .2ex}{\large\bf}}
\makeatother

\font\eufm=eufm10\font\eufms=eufm10\font\eufmss=eufm10\newfam\eufam
\textfont\eufam=\eufm\scriptfont\eufam=\eufms\scriptscriptfont\eufam=\eufmss

\input{epsf}

\begin{document}
\thispagestyle{empty}
\title{\bf Stabilization of a piezoelectric system}
\author{Ka\"{\i}s Ammari
\thanks{D\'epartement de Math\'ematiques, Facult\'e des Sciences de
Monastir, 5019 Monastir, Tunisie, \ email: kais.ammari@fsm.rnu.tn}
\, and \, Serge Nicaise
\thanks{Universit\'e de Valenciennes et du Hainaut Cambr\'esis,
LAMAV, FR CNRS 2956, Le Mont Houy, 59313 Valenciennes Cedex 9,
France, \, email : snicaise@univ-valenciennes.fr }}
\date{}
\maketitle
%
%
%============================INTRODUCTION============================%%
%
{\bf Abstract.~} {\small We consider a stabilization problem for a
piezoelectric system. We prove an exponential stability result under
some Lions geometric condition. Our method is based on an identity with multipliers
 that allows to show an appropriate observability estimate.\\

 }
\section{Introduction} \label{formulare}

We consider the dynamical  behavior of a piezoelectric system (which
means the ability of some materials, like ceramics and quartz, to
generate an electric field in response to applied mechanical
stress), where a proper modeling involves the displacement vector,
the electric field and the magnetic field, which are governed by the
elasticity system coupled with  Maxwell's equations. This system
plays an important role in various applications in structural
mechanics and in mechatronics, for such a model we refer to
\cite{EringenMaugin,Ikeda}.

Let $\Omega $ be a bounded domain of ${\mathbb R}^3$ with a
Lipschitz boundary $\Gamma $. In that domain we   consider the
non-stationary piezoelectric   system that consists in a coupling
between the elasticity system   with the Maxwell equation. More
precisely we analyze the  partial differential equations based on
the following relations between the stress tensor, the electric
displacement  and the magnetic induction: \be \label{eq1}
\sigma_{ij}(u,E) = a_{ijkl}\gamma_{kl}(u) - e_{kij}E_k \ \forall \,
i,j=1,2,3, \ee \be \label{eq2} D_i = \varepsilon_{ij}E_j + e_{ikl}
\gamma_{kl} (u) \, \forall \, i=1,2,3
 \ee
 \be
\label{eq3} B = \mu H.
\ee

The equations of equilibrium are \be \label{eq4bis}
\partial_t^2 u_i = \partial_j \sigma_{ji} \, \forall \, i=1,2,3
\ee
for the elastic displacement and
\be
\label{elastic}
\partial_t D = \rot  H, \, \partial_t B = - \rot E
\ee for the electric/magnetic fields.

This system models the coupling between Maxwell's system and the
elastic one, in which $E(x,t), H(x,t)$ are the electric and magnetic
fields at the point $x \in \Omega$ at time $t,\, u(x,t)$ is the
displacement field at the point $x \in \Omega$ at time $t,$ and
$\gamma_{ij}(u)^3_{i,j=1}$ is the strain tensor given by
$$
\gamma_{ij}(u) = \frac{1}{2} \, \left( \frac{\partial u_i}{\partial
x_j} + \frac{\partial u_j}{\partial x_i} \right).
$$
Here $\sigma=(\sigma_{ij})^3_{i,j=1},D=(D_1, D_2, D_3),$ and
$B=(B_1, B_2, B_3)$ are the stress tensor, electric displacement,
and magnetic induction, respectively. $\varepsilon,\mu$ are the
electric permittivity and magnetic permeability, respectively, and
we will assume that they are positive real numbers. The elasticity
tensor $(a_{ijkl})_{i,j,k,l=1,2,3}$ is made of constant entries such
that
$$
a_{ijkl} = a_{jikl} = a_{klij}
$$
and satisfies the ellipticity condition \be\label{ellip}
a_{ijkl}\gamma_{ij}\gamma_{kl} \geq \alpha_0 \gamma_{ij}
\gamma_{ij}, \ee for every symmetric tensor $(\gamma_{ij})$ and some
$\alpha_0 > 0$. The piezoelectric tensor $e_{kij}$ is also made of
constant entries such that
$$
e_{kij} = e_{kji}.
$$
For shortness in the remainder of the paper introduce the tensor
$\sigma(u)=(a_{ijkl}\gamma_{kl}(u))^3_{i,j=1}$ and let
 $\nabla \sigma$
be the vector field defined by
$$
\nabla \sigma =(\partial_j \sigma_{ij})_{i=1}^3,
$$
while for a tensor $\gamma=(\gamma_{ij})^3_{i,j=1},$ and a vector
$F=(F_1, F_2, F_3),$ we set
$$
e\gamma=(e_{ikl}\gamma_{kl})^3_{i=1} \quad e^\top F=
(e_{ikl}F_{i})^3_{k,l=1}.
$$
These last notations mean that  $e$ corresponds to a linear mapping
from $\R^{3\times 3}$ into $\R^{3}$ and that $e^\top$ is its
adjoint. With these notations, we see that (\ref{eq1}) is equivalent
to
$$
\sigma(u,E)=\sigma(u)-e^\top E,
$$
while (\ref{eq2}) is equivalent to
$$
D = \varepsilon E  + e  \gamma(u).
$$

The system (\ref{eq1})-(\ref{eq3}) is completed with the boundary
and Cauchy conditions. This means that we are considering the
following system \be\label{max} \left\{
\begin{tabular}{llllllll}
& $\partial_t^2 u-\nabla \sigma(u,E) =0 $ in $Q:=\Omega \times
]0,+\infty[$,\vspace{2mm}
\\
& ${\partial_t D} -\rot H=0 $ in $Q$,\vspace{2mm}
\\
& $\mu  \partial_t H +\rot E=0 $ in $Q$,\vspace{2mm}
\\
& $\ddiv(D)=\ddiv(\mu H)=0$ in $Q$,\vspace{2mm}
\\
 &$H
\times \nu - (Q^*\partial_t u) \times \nu + (E \times \nu) \times
\nu = 0$  on $\Sigma= \Gamma \times (0,+\infty),$ \vspace{2mm}
\\
&$\sigma(u,E)\cdot \nu  + Q (E\times \nu)+Au+
\partial_t u = 0$ on $\Sigma,$\vspace{2mm}
\\
& $u(0) =  u_0,$ $\partial_t u(0) =  u_1$ in $\Omega $,\vspace{2mm}
\\
& $E(0) =  E_0,$ $H(0) =  H_0$ in $\Omega $,
\end{tabular}
\right. \ee   where $\nu$ is the unit normal vector of $\partial
\Omega$ pointing towards the exterior of $\Omega$, $A$ is a positive
constant and $Q$ is a function from $\Gamma$ into the set of
$3\times 3$ matrices with the regularity $Q\in L^\infty(\Gamma,
\C^{3\times 3})$.

\br\label{rQnormal} {\rm Note that the image of $Q$ of normal vector
fields plays no role in the boundary conditions appearing in
(\ref{max}). Indeed for $X\in \C^3$, let $X_\nu=(X\cdot \nu) \nu$
and $X_\tau=X-X_\nu$ be the normal and tangential components of $X$
respectively, by writing
\[
Q_\nu X=Q X_\nu, \quad Q_\tau X=Q X_\tau,
\]
we get the splitting
\[
Q=Q_\nu+Q_\tau, \quad Q^*=Q_\nu^*+Q_\tau^*.
\]
But by definition $Q_\nu^*X$ is orthogonal to the tangent plane.
Therefore $Q (E\times \nu)=Q_T (E\times \nu)$  and $(Q^*\partial_t
u) \times \nu=(Q^*_\tau\partial_t u) \times \nu$, which means that
the normal part $Q_\nu$ of $Q$ does not contribute to the boundary
conditions.} \er

Boundary or internal stability of the second order elliptic systems,
like the wave equation or the elasticity system, have been studied
by many authors, let us quote
\cite{Chen,Alabau,Bey,Guesmia,Horn,Komornikbook,LaJDE83,LaSiam83,LaTri,Martinez}
among others. Similar results for Maxwell's system can be found in
\cite{Kapi,Komo,Komo98,Komornikbook,ELN,ELN2,La89,Phung}. The
combination of these results to the piezoelectric system in some
particular cases has been treated in
\cite{KapiRaupp,kapi03,nicaise}. For the quasi-static case
(corresponding to the hypothesis that $E$ is curl free, hence the
gradient of a potentiel), we can refer to
\cite{kapitonov:06,miara:09,miara:09b}.

%ici

In \cite{kapi03} the authors consider the above problem in the case
$Q=0$ with eventually discontinuous coefficients  and an additional
memory term and prove the exponential decay rate of the energy if
$A$ is small enough and if $\Omega$ satisfies some geometrical
conditions (star like shape). On the contrary in \cite{nicaise}, the
author treats the case $Q=I$ and some nonlinear feedback terms, but
with the choice of $e$ such that $\nabla (e^\top E)=\xi \rot E$ for
some real number $\xi$ (case excluding the natural condition
$e_{kij} = e_{kji}$) and proves  the exponential decay rate of the
energy in the case of linear feedbacks if $\Omega$ is strictly star
shaped with respect to a point. In that last paper the author
combines the multiplier technique with the one from \cite{Bey},
where the authors uses some tangential integration by parts and a
technique from \cite{conrad:93}. Our goal is here to perform the
same analysis for the general system (\ref{max}). For $Q\in
L^\infty(\Gamma, \C^{3\times 3})$  we   prove that the system is
well-posed using semigroup theory. On the other hand using the
multiplier method (see \cite{kapi03}) and a technique inspired from
\cite{Alabau,conrad:93,HNicSene,nicaise} to absorb a zero order
boundary term, we show that the system is exponential stable if
$Q=\alpha \, I$ for some scalar continuously differentiable function
$\alpha$ such that $\nabla \alpha$ is small enough.

The paper is organized as follows. The second section deals with the
well-posedness of the problem. In the last section we give the main
result of this paper which is the exponential stability of the
piezoelectric system and its proof.

\section{Well-posedness of the problem\label{well}}

We start this section with the well-posedness of  problem
(\ref{max}). At the end we will check the dissipativeness of
(\ref{max}).

Let us  introduce the Hilbert spaces (see e.g. \cite{La89,Nicmax})
\beq\label{2.0} && J(\Omega)=\{E \in L^2(\Omega )^3\vert \ddiv E=0
\hbox{ in } \Om\},
\\
&&\caH=H^1(\Om)^3\times L^2(\Om)^3\times L^2(\Om)^3\times J(\Omega),
\eeq equipped with the  norm induced by the inner product \beqs &&
(E,E')_\eps=\int_\Om\eps \, E(x)\cdot E'(x)\,dx, \forall E,E'\in
J(\Omega ),
\\
&&
\left((u,v,E,H),(u',v',E',H')\right)_\caH=(u,u')_1+(v,v')_0
\\
&&+(E,E')_\eps+(H,H')_\mu, \forall  (u,v,E,H),(u',v',E',H')\in \caH,
\eeqs where we have set \beqs &&(u,u')_0=\int_\Om u(x)\cdot
u'(x)\,dx,
\\
&&(u,u')_1=\int_\Om \sigma(u)(x): \gamma(u')(x)\,dx+A\int_\Gamma
u(x)\cdot u'(x)\,dS, \eeqs with the notation
$$
\sigma(v):\gamma(v'):= \sigma_{ij}(v)\gamma_{ij}(v').
$$

Now define the linear operator ${\mathcal A}$ from $\caH$ into
itself as follows: \beq\label{domA} D({\mathcal
A})&=&\{\left(u,v,E,H \right)\in \caH\vert \nabla \sigma(u,E), \rot
E, \rot H\in L^2(\Omega )^3; v\in H^1(\Om)^3;
\\
\nonumber
&&E\times \nu, H\times \nu \in L^2(\Gamma )^3 \hbox{ satisfying }
\\
\label{bc1} &&H\times \nu -(Q^*v) \times \nu +E\times \nu \times \nu
= 0 \hbox{ on } \Gamma ,
\\
\label{bc2} &&\sigma(u,E) \cdot\nu+Au +v+Q (E\times \nu)  = 0 \hbox{
on } \Gamma \}. \eeq    For all $(u,v,E,H) \in D({\mathcal A})$ we
take
$$
{\mathcal A}(u,v,E,H)= \left(v, \nabla\sigma(u,E), \varepsilon
^{-1}(\rot H-e\gamma(v)),-\mu ^{-1}\rot E \right).
$$

The boundary conditions (\ref{bc1})  and (\ref{bc2}) are meaningful
since  for $(u,v,E,H) \in D({\mathcal A})$,
 from section 2 of \cite{ABDG}
the property $\nabla\sigma(u,E)\in L^2(\Om)^3$ implies that $
\sigma(u,E)\cdot \nu$ belongs to $H^{-1/2}(\Gamma)^3$. Since the
properties $u,v \in H^1(\Om)^3$ imply that $Au+v$ belongs to
$H^{1/2}(\Gamma)^3$, the boundary condition (\ref{bc2}) has a
meaning (in  $H^{-1/2}(\Gamma)^3$) and furthermore yields
$\sigma(u,E)\cdot \nu\in L^{2}(\Gamma)^3$ (because $Q (E\times
\nu)\in L^{2}(\Gamma)^3$). Similarly the properties of $H$ and $v$
give a meaning to the boundary condition (\ref{bc1}) (as an equality
in $L^{2}(\Gamma)^3$). In summary both   boundary conditions
(\ref{bc1}) and (\ref{bc2}) have to be understood as an equality in
$L^{2}(\Gamma)^3$.

We now see that formally problem (\ref{max}) is equivalent to
\be\label{3.1bis}
\left\{
\begin{tabular}{ll}
& $\frac{\partial U}{\partial t}={\mathcal A}U,$ \\
& $U(0) = U_0,$
\end{tabular}
\right.
\ee
when $U =(u,\partial_t u,E,H)$ and $U_0=(u_0,u_1,E_0,H_0)$.

We shall  prove that this problem (\ref{3.1bis}) has a unique
solution using  semigroup theory   by showing that ${\mathcal A}$ is
a maximal dissipative operator.

\bl\label{lmaxmon}
 ${\mathcal A}$ is a maximal dissipative operator.
\el \pf We start with the dissipativeness:
$$
({\mathcal A}U,U)_\caH\leq 0, \forall \, U\in D({\mathcal A}).
$$
From the definition of ${\mathcal A}$ and the inner product in
$\caH$,  we have \beqs &&({\mathcal
A}U,U)_\caH=(v,u)_1+(\nabla\sigma(u,E), v)_0
\\
&&+ \int_\Om \{E\cdot(\rot H-e\gamma(v))  - \rot E\cdot H\} \,dx,
\eeqs for any $(u,v,E,H)\in D({\mathcal A})$. Lemma 2.2 of
\cite{ELN2} and Green's formula yield equivalently \beqs
&&({\mathcal A}U,U)_\caH=(v,u)_1- \int_\Om\sigma(u,E):\gamma( v)\,dx
\\
&&-\int_\Om e\gamma(v)\cdot E \,dx
\\
&&+\int_{\Gamma} \{(\sigma(u,E)\cdot \nu)\cdot v + (E\times
\nu)\cdot H\}\,dS, \eeqs for any $(u,v,E,H)\in D({\mathcal A})$.
Using the definition of the inner product $(\cdot, \cdot)_1$ and the
boundary conditions (\ref{bc1}) and (\ref{bc2}), we arrive at
$$
({\mathcal A}U,U)_\caH=-\int_\Gamma\{|v|^2+|E\times \nu|^2\}
\,dS\leq 0,
$$
for any $(u,v,E,H)\in D({\mathcal A})$.

Let us now pass to the maximality. This means that for at least one
non negative real number $\lambda$, $\lambda \,I-{\mathcal A}$ has
to be surjective. Let us show that indeed  $I-{\mathcal A}$  is
surjective. This means that for all $\left(f,g,F,G \right)$ in
$\caH$, we are looking for $\left(u,v,E,H \right)$ in $D({\mathcal
A})$ such that \be\label{surj} (I-{\mathcal A}) \left(u,v,E,H
\right)= \left(f,g,F,G \right). \ee From the definition of
${\mathcal A}$, this equivalently means \be\label{surjbis} \left\{
\begin{array}{llll}
&&u-v=f,\\
&&v-\nabla\sigma(u,E)=g,\\
&&E-\eps^{-1}(\rot H-e\gamma(v))=F,\\
&&H+\mu ^{-1}\rot E=G.
\end{array}
\right.
\ee
The first and fourth equations allow to eliminate $H$ and $v$, since
they are respectively equivalent to \beq\label{eliv}
&&v=u-f,\\
\label{psi} &&H =G-\mu ^{-1}\rot E. \eeq Substituting these
expressions in the second and third equations yields formally \beq
\label{equ}
&& u-\nabla\sigma(u,E)=f+g,\\
\label{eqE} &&\eps E +\rot(\mu ^{-1}\rot E )+e\gamma(u) =\eps F+\rot
G+e\gamma(f). \eeq

This   system in $(u,E)$ will be uniquely defined by  adding
boundary conditions on $u$ and $E$. Indeed using
 the identities  (\ref{eliv}) and (\ref{psi}), we see that (\ref{bc1})
and (\ref{bc2}) are formally equivalent to \beq\label{bconE} &&-\mu
^{-1}\rot E \times \nu+Q^*  u\times \nu+(E \times \nu) \times \nu=-G
\times \nu+Q^* f \times \nu \hbox{ on } \Gamma,
\\
\label{bconv} &&\sigma(u,E)\cdot\nu+Au+u+Q(E\times \nu)=f \hbox{ on
} \Gamma. \eeq

By formal integration by parts we  remark that the variational formulation of the system (\ref{equ})-(\ref{eqE})
 with the boundary conditions (\ref{bconE})-(\ref{bconv}) is the following one:
Find $(u,E)\in V$ such that \be\label{varmax}
a((u,E),(u',E'))=F(u',E'), \forall (u',E') \in V, \ee where the
Hilbert space $V$ is given by $V=H^1(\Om)^3\times W$ when
$W$ is defined by
$$
W=\{E\in L^2(\Omega )^3\vert \rot E\in L^2(\Omega )^3  \hbox{
and } E\times \nu \in L^2(\Gamma )^3\},
$$
with the norm
$$
|| E||_{W}^2=\int_\Om (|E|^2+|\rot E|^2)dx+\int_\Gamma |E\times \nu
|^2\,dS,
$$
 the form $a$ is defined by
\beqs a((u,E),(u',E'))&=& \int_\Omega
\{\sigma(u,E):\gamma(u')+u\cdot u'\}\;dx
\\
&+& \int_\Omega \{\mu ^{-1}\rot E\cdot\rot E' + \eps E\cdot E'+
e\gamma(u)\cdot E'\}\;dx
\\
&+&\int_\Gamma \{(E\times\nu)\cdot (E'\times \nu)+(A+1)u\cdot
u'+Q(E\times \nu)\cdot u'-Q^*u\cdot (E'\times \nu)\}\,dS, \eeqs and
finally  the form $F$ is defined by
$$
F(u',E')= \int_\Omega \{(f+g)\cdot u'+ (\eps F+e\gamma(f))\cdot E' +
G\cdot\rot E'\}\;dx +\int_\Gamma (f\cdot u'-(Q^*f\times \nu)\cdot
E')\, dS.
$$

We easily see that the bilinear form $a$ is coercive on $V$ since
\beqs a((u,E),(u,E))&=& \int_\Omega
\{\sigma(u):\gamma(u)+|u|^2\}\;dx
\\
&+& \int_\Omega \{\mu ^{-1}|\rot E|^2+ \eps |E|^2\}\;dx
\\
&+&\int_\Gamma \{|E\times\nu|^2+(A+1)|u|^2\}\,dS, \eeqs which is
clearly greater than $||u||^2_{H^1(\Omega)^3}+|| E||_{W}^2$ by the
ellipticity assumption on the elasticity tensor. Hence by the
Lax-Milgram lemma, problem (\ref{varmax}) has a unique solution
$(u,E)\in V$.

To end our proof we need to show that the solution $(u,E)\in V$
 of (\ref{varmax})  and $v$, $H$ given respectively  by (\ref{eliv}),
(\ref{psi}) are such that
  $(u,v,E,H)$ belongs to  $D({\mathcal A})$  and satisfies (\ref{surj}) (or equivalently (\ref{surjbis})).
First taking  test functions $u'$ in $\caD(\Om)^3$ and $E'=0$, we
get
$$
\nabla \sigma(u,E)+v=g \hbox{ in } \caD'(\Om).
$$
This implies the second identity in (\ref{surjbis}) as well as the
regularity $ \nabla \sigma(u,E)\in L^2(\Om)^3$ (from the fact that
$v, \rot E$ as well as $g$ belongs to that space).

Second we take test functions $u'=0$ and $E'=\chi$ with
$\chi\in\caD(\Om)^3$ by Lemma 2.3 of \cite{ELN2} we get
$$
\eps E-\rot H+e\gamma(u)=\eps F \hbox{ in } \caD'(\Om).
$$
This means that the third identity in (\ref{surjbis}) holds as well
as the regularity $\rot H\in L^2(\Om)^3$.

Thirdly taking test functions $v'\in H^1(\Om)^3$ and $E'=\chi$ with
$\chi\in C^\infty(\bar\Om)^3$  and applying Green's formula (see
section 2 of \cite{ABDG} and Lemma 2.2 of \cite{ELN2}), we get \beqs
&&\langle \sigma(u,E)\cdot \nu,v'\rangle-\int_\Gamma (H\times
\nu)\cdot E'\,dS +\int_\Gamma (Q(E\times \nu)\cdot u'-(Q^*u\times
\nu)\cdot E'\,dS
\\
&& +\int_\Gamma \{(E\times\nu)\cdot (E'\times \nu)+(A+1)u\cdot
u'\}\,dS =0. \eeqs This leads to the boundary conditions (\ref{bc1})
and (\ref{bc2}) since $u'$ (resp. $\chi$) was arbitrary in
$H^1(\Om)^3$  (resp. in $C^\infty(\bar\Om)^3$) whose   trace belongs
to a dense subspace of  $L^{2}(\Gamma)^3$.

Finally from (\ref{psi}) and the fact that $\mu G$ is divergence free, $\mu H$ is also divergence free.
\findem

Semigroup theory
 \cite{Pazy,Showalter} allows to conclude the following existence results:
\bc\label{c1.4} For all $(u_0,u_1,E_0,H_0)\in \caH$, the problem
(\ref{max}) admits a unique (weak) solution $\left(u,E,H\right)$
satisfying $(u,\partial_tu,E,H)\in C(\R_+,\caH)$, or equivalently
$u\in C^1(\R_+,L^2(\Om)^3)\cap C(\R_+,H^1(\Om)^3)$, $E\in
C(\R_+,L^2(\Om)^3)$ and $H\in C(\R_+,J(\Om))$. If moreover
$\left(u_0,u_1,E_0,H_0 \right)$ belongs to $D({\mathcal A})$ and
satisfies $$ \ddiv (e\gamma(u_0)+\eps E_0)=0 \hbox{ in } \Om,
$$
then the problem (\ref{max}) admits a unique (strong) solution
$\left(u,E,H\right)$ satisfying $(u,\partial_tu,E,H)\in
C^{1}(\R_+,\caH)$ $\cap C(\R_+,D({\mathcal A}))$, or equivalently
satisfying $u\in C^{2}(\R_+,L^2(\Om)^3)\cap C^{1}(\R_+,H^1(\Om)^3)$,
$E\in C^{1}(\R_+,J(\Om))\cap C(\R_+,W)$,
 $H\in C^{1}(\R_+,J(\Om))\cap
C(\R_+,W)$,   satisfying (\ref{bc1})-(\ref{bc2}) for a.e. $t$ (with
$v=\partial_tu$), as well as
$$
 \nabla\sigma(u,E)\in C(\R_+,L^2(\Om)^3).
$$
Note that, in that last case,   $D=e\gamma(u)+\eps E$  satisfies in
particular
$$ \ddiv D=0 \hbox{ in } \Om\times \R_+.
$$
\ec

We finish this section by showing the dissipativeness of our system.
\bl\label{l2.1} The energy \beq\label{2.1} {\mathcal
E}(t)&=&\frac{1}{2}\int_\Om(|\partial_tu(x,t)|^2+\sigma(u)(x,t):\gamma(u)(x,t))\,dx
+\frac{A}{2}\int_\Gamma|u(x,t)|^2\,dS(x)
\\
\nonumber &+& \frac{1}{2}\int_\Om(\eps |{\mathcal E}(x,t)|^2+\mu
|H(x,t)|^2)\,dx \eeq is non-increasing. Moreover
 for $\left(u_0,u_1,
E_0, H_0 \right)\in D({\mathcal A})$, we have
 for all $0\leq S<T<\infty$
\be\label{2.2weak} {\mathcal E}(S)-{\mathcal
E}(T)=\int_S^T\int_\Gamma \{|E(x,t)\times \nu|^2+|\partial_t
u(x,t)|^2\} \,dS dt, \ee and for all $t\geq 0$ \be\label{2.2}
\partial_t{\mathcal E}(t)=-\int_\Gamma \{|E(x,t)\times
\nu|^2+|\partial_t u (x,t)|^2\}\,dS. \ee \el \pf Since $D({\mathcal
A})$ is dense in $\caH$ it suffices to show (\ref{2.2}). For
$\left(u_0,u_1, E_0, H_0 \right)\in D({\mathcal A})$, from the
regularity of $u,E,H$, we have \beqs
\partial_t {\mathcal E}(t)&=&\int_\Om\{\partial_t^2u \cdot
\partial_tu +\sigma(u):\gamma(\partial_tu)\}\,dx +A\int_\Gamma
\partial_tu \cdot u \,dS
\\
&+& \int_\Om\{\eps E \cdot \partial_t E +\mu H \cdot
\partial_tH\}\,dx. \eeqs By (\ref{max}), we get \beqs
\partial_t {\mathcal E}(t)&=&\int_\Om\{\partial_tu \cdot \nabla\sigma(u,E)
+\sigma(u):\gamma(\partial_tu)\}\,dx +A\int_\Gamma \partial_t u
\cdot u \,dS
\\
&+&\int_\Om\{E\cdot (\rot H-E\gamma(\partial_tu)\}\,dx
\\
&=&\left(A\left( u(t), \partial_t u(t),E(t), H(t) \right), \left(
u(t), \partial_t u(t),E(t), H(t) \right)\right)_\caH. \eeqs We
conclude by Lemma \ref{lmaxmon}. \findem

\section{Exponential stability\label{expstab}}
In this section we prove the main result of this paper, namely the
exponential stability of our system (\ref{max}) when $\Om$ is
strictly star-shaped with respect to a point $x_0$. This result is
based on an identity with multipliers proved in \cite{kapi03} that
allows to show the next observability estimate.

\bt\label{tobsest} Assume that there exists $x_0\in\R^n$ and
$\delta>0$ such that \be\label{starshaped}
 m(x)\cdot \nu(x)\geq \delta \quad \forall x\in
\partial\Om,
\ee where $m(x)=x-x_0$. Assume also that $Q=\alpha \,I$ with a
continuously differentiable  function $\alpha$ from $\Gamma$ to
$\C$. Set $c_\alpha= \ds \max_\Omega|\nabla \alpha|$. Let
$\left(u,E,H\right)$ be the strong solution  of problem (\ref{max}).
Then there exists a positive constants $C$ (independent of $\alpha$)
such that for all $T>0$, and all $\theta$, there exists a constant
$C(\theta)$ (independent of $T$) such that the next observability
estimate holds: \be\label{obsest} T {\mathcal E}(T)\leq
\big(C(\theta) (1+c_\alpha T)+\theta T\big) {\mathcal E}(0)+ C
\int_{\Sigma_T} (|\partial_t u|^2+|E_\tau|^2)\,dSdt,\ee where
$\Sigma_T= \Gamma \times (0,T)$.
 \et \pf
 First the identity (3.9) of \cite{kapi03} with $t_0=0$ and
 $\varphi(x)=|x-x_0|^2/2$ yields
\be\label{3.9kapi}
 T {\mathcal E}(T)=r+\int_{\Sigma_T} V(x,t)\,dS(x)dt,
\ee where we have set  \beqs r&=&-2\int_\Om \Big\{\partial_t u\cdot
\{u+(m\cdot \nabla) u)\}+\mu (m\times H)\cdot \{ \eps E
+e\gamma(u)\}\Big\}\,dx\Big|_0^T,
\\
V&=&2 \{t \partial_t u+(m\cdot \nabla) u+u\}\cdot \sigma(u,E)\nu
+m\cdot\nu \{|\partial_tu|^2-\sigma(u):\gamma(u)+\eps |E|^2+\mu
|H|^2\}
\\
&+&2t (H\times E)\cdot \nu-2\eps E\cdot \nu E\cdot m-2\mu H\cdot \nu
H\cdot m
\\
&-&2(m\times e \gamma(u))\cdot (E\times \nu). \eeqs Using the
boundary conditions from (\ref{max}), we see that \beqs V&=&-2
t\partial_tu(Q(E\times \nu)+Au+\partial_tu) +\Delta
\\
&+&m\cdot\nu \{\eps |E|^2+\mu |H|^2\}
\\
&-&2t (Q^*\partial_t u) \times \nu)\cdot E-2t |E_\tau|^2-2\eps E_\nu
(E_\nu m\cdot \nu+E_\tau \cdot m_\tau)-2\mu H_\nu (H_\nu m\cdot
\nu+H_\tau \cdot m_\tau)
\\
&-&2(m\times e \gamma(u))\cdot (E\times \nu), \eeqs where we recall
that $E_\nu=E\cdot \nu$, $E_\tau=E-E_\nu\nu$ and
$$
\Delta=2\{(m\cdot \nabla) u+u\}\cdot \sigma(u,E)\nu +m\cdot\nu
\{|\partial_tu|^2-\sigma(u):\gamma(u)\}.
$$
By Young's inequality, there exists $C>0$ such that for all $
\beta_1, \beta_2>0$ \beqs V&\leq& -2 At u \partial_tu
  \\
&-&(m\cdot \nu-\beta_2)(\eps |E_\nu|^2+\mu |H_\nu|^2)
  \\
  &+&\Delta+ C(1+\frac{1}{\beta_2}+\frac{1}{\beta_1})|E\times \nu|^2+C(1+\frac{1}{\beta_2})|H\times
  \nu|^2+\beta_1\gamma(u):\gamma(u).
 \eeqs
 By using    again the
 first boundary condition from (\ref{max}), we get for all
 $\beta_1,\beta_2>0$
\beq\label{estV} V&\leq& -2 At u \partial_tu+C(1+\frac{1}{\beta_2})
|\partial_tu|^2+C(1+\frac{1}{\beta_2}+\frac{1}{\beta_1})|E\times
\nu|^2
 \\
&-&(m\cdot \nu-\beta_2)(\eps |E_\nu|^2+\mu |H_\nu|^2)
  +\Delta+ \beta_1\gamma(u):\gamma(u).
  \nonumber
 \eeq

 Let us   transform the first term of this right-hand side:
 $$
 -2 A\int_{\Sigma_T}  t u \partial_tu\, dSdt =- A\int_{\Sigma_T}  t \frac{d}{dt}
 u^2\, dSdt,
 $$
 and by an integration by parts in time, we get
 $$
 -2 A\int_{\Sigma_T}  t u \partial_tu\, dSdt=  A\int_{\Sigma_T}
 u^2\, dSdt-A\int_{\Omega}  t
 u^2\,dx\Big|_0^T.
 $$
 This proves that
 \be\label{estt1}
-2 A\int_{\Sigma_T}  t u \partial_tu\, dSdt\leq  A\int_{\Sigma_T}
 u^2\, dSdt.
 \ee
Let us now estimate the term $\Delta$. First using  the
 second boundary condition from (\ref{max}), we see that
$$ \Delta=-2\{(m\cdot \nabla) u+u\}\cdot (Q (E\times \nu)+Au+\partial_tu)
+m\cdot\nu \{|\partial_tu|^2-\sigma(u):\gamma(u)\}. $$ Using the
ellipticity assumption (\ref{ellip}) and  condition
(\ref{starshaped}) we obtain \beq\label{estdelta} \quad \quad
\quad\Delta&\leq& -2\{(m\cdot \nabla) u+u\}\cdot (Q (E\times
\nu)+Au+\partial_tu) +m\cdot\nu |\partial_tu|^2-\alpha_0 \delta
\gamma(u):\gamma(u) \\
&\leq &-2u\cdot Q (E\times \nu)-2A |u|^2-2 u\cdot \partial_tu
-2(m\cdot \nabla) u\cdot Q (E\times \nu)
 \nonumber
\\
 &-&2A(m\cdot \nabla) u\cdot  u
 -2 (m\cdot \nabla) u \cdot  \partial_tu +m\cdot\nu |\partial_tu|^2-\alpha_0 \delta
\gamma(u):\gamma(u).\nonumber \eeq We need to estimate some terms of
this right-hand side. First as before an integration by parts in
time yields
$$
-2 \int_{\Sigma_T}   u \partial_tu\, dSdt\leq  A\int_{\Gamma}
 |u(x, t=0)|^2\,dS(x)\leq 2 {\mathcal E}(0).
$$

As in \cite{Bey,HNicSene}, one can show that \beq\label{sn1} -2
A\int_{\Sigma_T} (m\cdot \nabla) u\cdot u \,dS dt &\leq&
\frac{C}{\theta_1} \int_{\Sigma_T} |u|^2\, dSdt+ \theta_1
\int_{\Sigma_T} \gamma(u):\gamma(u)\,dS dt, \eeq
  as well as
\be\label{sn1bis} \int_{\Sigma_T} (m\cdot \nabla) u\cdot
\partial_tu\,dS dt \leq C {\mathcal E}(0)+ \frac{C}{\theta_2}
\int_{\Sigma_T} (|u|^2+|\partial_tu|^2)\, dSdt+ \theta_2
\int_{\Sigma_T} \gamma(u):\gamma(u)\,dS dt, \forall
\theta_1,\theta_2>0.\ee

 By Young's
inequality we clearly have \be\label{sn2} \int_{\Sigma_T} u\cdot
Q(E\times \nu)\, dSdt \leq C \int_{\Sigma_T} (|u|^2+|E_\tau|^2)\,
dSdt. \ee
 Now we notice that
$$
(m\cdot \nabla) u\cdot Q(E\times \nu)=(Q^*(m\cdot \nabla)
u)\cdot(E\times \nu),
$$
and for any $k=1,2,3$, we may write \beqs (Q^*(m\cdot \nabla)
u)_k&=&Q^*_{kj} m_i \partial_iu_j \\
&=& 2Q^*_{kj} m_i \gamma_{ij}(u)-Q^*_{kj} m_i \partial_ju_i
\\
&=& 2Q^*_{kj} m_i \gamma_{ij}(u)  +Q^*_{kj} u_i
\partial_jm_i-Q^*_{kj} \partial_j (m_i  u_i).
\eeqs The two first terms of this right-hand side will be estimated
by Young's inequality and it therefore remains to estimate the last
term, namely by the previous identities we have \beq\label{sn3}
&&\int_{\Sigma_T} (m\cdot \nabla) u\cdot Q(E\times \nu) \, dSdt \leq
\theta_3 \int_{\Sigma_T} \gamma(u):\gamma(u)\,dS dt \\
&& + \, C  \int_{\Sigma_T}
(|u|^2+(1+\frac{1}{\theta_3})|E_\tau|^2)\, dSdt -\int_{\Sigma_T}
(Q^* \nabla (m\cdot u))\cdot (E\times \nu)\, dSdt, \forall
\theta_3>0. \nonumber \eeq Now using Green's formula, we see that
$$\int_{\Sigma_T} (Q^* \nabla (m\cdot u))\cdot (E\times \nu)\, dSdt =
\int_{Q_T} \{\rot(Q^* \nabla (m\cdot u))\cdot E-Q^* \nabla (m\cdot
u) \cdot \rot E\}\, dSdt,
$$
where $Q_T=\Omega\times (0,T)$. Now using the fact that
$Q(x)=\alpha(x) \,I,$ and that $\rot E=\mu
\partial_t H$, we obtain
$$\int_{\Sigma_T} (Q^* \nabla (m\cdot
u))\cdot (E\times \nu)\, dSdt = \int_{Q_T} \{(\nabla \alpha \times
\nabla (m\cdot u))\cdot E-Q^* \nabla (m\cdot u) \cdot \mu \partial_t
H\}\, dSdt.
$$
For this last term, we first integrate by parts in time and get
$$
\int_{Q_T} Q^* \nabla (m\cdot u) \cdot \mu \partial_t H\, dSdt=-
\int_{Q_T} Q^* \nabla (m\cdot \partial_tu) \cdot \mu H\,
dSdt+\int_{\Omega} Q^* \nabla (m\cdot u) \cdot \mu H\,dx\Big|_0^T.
$$
An integration by parts in space leads to \beqs \int_{Q_T} Q^*
\nabla (m\cdot u) \cdot \mu \partial_t H\, dSdt&=& \int_{Q_T} (Q^*
 m\cdot \partial_tu \ddiv( \mu H)+m\cdot \partial_tu \nabla \alpha\cdot(\mu H))\, dSdt
 \\
 &-&\int_{\Sigma_T} Q^*
 m\cdot \partial_tu   (\mu H)\cdot\nu \, dSdt
 +\int_{\Omega} Q^* \nabla (m\cdot u) \cdot
\mu H\,dx\Big|_0^T. \eeqs These two identities and reminding that
$\ddiv( \mu H)=0$ lead to \beqs \int_{\Sigma_T} (Q^* \nabla (m\cdot
u))\cdot (E\times \nu)\, dSdt &=& \int_{Q_T} \{(\nabla \alpha \times
\nabla (m\cdot u))\cdot E- m\cdot \partial_tu \nabla \alpha\cdot(\mu
H)\}\, dSdt
\\
&+&\int_{\Sigma_T} Q^*
 m\cdot \partial_tu   (\mu H)\cdot\nu\, dSdt
 -\int_{\Omega} Q^* \nabla (m\cdot u) \cdot
\mu H\,dx\Big|_0^T. \eeqs By Young's inequality we find that \beqs
\int_{\Sigma_T} (Q^* \nabla (m\cdot u))\cdot (E\times \nu) \,
dSdt\leq C(1+c_\alpha T) {\mathcal E}(0)+  \int_{\Sigma_T}
\{\frac{C}{\theta_4}|\partial_tu|^2+\theta_4 |H_\nu|^2\}\, dSdt,
\forall \theta_4>0. \eeqs  This last estimate in (\ref{sn3}) leads
to
\beq\label{sn4} &&\int_{\Sigma_T} (m\cdot \nabla) u\cdot Q(E\times
\nu)\, dSdt \leq C(1+c_\alpha T) {\mathcal E}(0)+\theta_3
 \int_{\Sigma_T} \gamma(u):\gamma(u)\, dS dt \\
&&+ C \int_{\Sigma_T}
\{|u|^2+\frac{1}{\theta_4}|\partial_tu|^2+(1+\frac{1}{\theta_3})|E_\tau|^2)+\theta_4
|H_\nu|^2\}\, dSdt, \forall \theta_3, \theta_4>0. \nonumber \eeq Now
using again Young's inequality and the estimates (\ref{sn1}),
(\ref{sn1bis}), (\ref{sn2}) and (\ref{sn4}) into the identity
(\ref{estdelta}), we obtain that
\beqs \int_{\Sigma_T}\Delta\, dSdt&\leq& C(1+c_\alpha T) {\mathcal
E}(0)+(-\alpha_0\delta+\theta_1+\theta_2+\theta_3)
 \int_{\Sigma_T} \gamma(u):\gamma(u)\, dS dt \\
&&+ C \int_{\Sigma_T}
\{(1+\frac{1}{\theta_2}+\frac{1}{\theta_4})|\partial_tu|^2+(1+\frac{1}{\theta_3})|E_\tau|^2)+\theta_4
|H_\nu|^2\}\, dSdt \nonumber \\
&+&C\{(1+\frac{1}{\theta_1}+\frac{1}{\theta_2}\}\int_{\Sigma_T}
|u|^2\, dSdt, \forall \, \theta_1, \theta_2, \theta_3, \theta_4 >0,
\nonumber
 \eeqs
This estimate in (\ref{estV}) and using (\ref{estt1}), we get
finally
\beqs \int_{\Sigma_T}V\, dSdt&\leq& C(1+c_\alpha T) {\mathcal
E}(0)+(-\alpha_0\delta+\beta_1+\theta_1+\theta_2 +\theta_3)
 \int_{\Sigma_T} \gamma(u):\gamma(u)\, dS dt \\
&&+ \, C \int_{\Sigma_T} (1+\frac{1}{\beta_2}+\frac{1}{\theta_2}
+\frac{1}{\theta_1}+\frac{1}{\theta_4})|\partial_tu|^2\, dSdt
\\
&&+ \, C \int_{\Sigma_T}
(1+\frac{1}{\beta_1}+\frac{1}{\beta_2}+\frac{1}{\theta_3})|E_\tau|^2\,
dSdt
\\
&&+ \, C \int_{\Sigma_T}\{(-m\cdot \nu+\beta_2)\eps |E_\nu|^2+
((-m\cdot \nu+\beta_2)\mu+\theta_4) |H_\nu|^2
\}\, dSdt \nonumber \\
&+&C \, \{1+\frac{1}{\theta_2}+\frac{1}{\theta_1}\}\int_{\Sigma_T}
|u|^2\, dSdt, \forall \beta_1,\beta_2,\theta_1,\theta_2,\theta_3,
\theta_4>0.\nonumber \eeqs

By choosing  $\beta_1$,  $\beta_2$, $\theta_1$, $\theta_2$,
$\theta_3$ and $\theta_4$ small enough,   we have found that
\be\label{estv} \int_{\Sigma_T}V\, dSdt\leq C (1+c_\alpha T)
{\mathcal E}(0)+ C \int_{\Sigma_T}
(|u|^2+|\partial_tu|^2+|E_\tau|^2)\, dSdt.\ee Coming back to
(\ref{3.9kapi}) and using again Young's and Korn's inequalities to
estimate $r$, we obtain \be\label{star} T {\mathcal E}(T)\leq C
(1+c_\alpha T) {\mathcal E}(0)+ C \int_{\Sigma_T}
(|u|^2+|\partial_tu|^2+|E_\tau|^2)\, dSdt.\ee Now invoking Lemma
\ref{l21} below, we arrive at \beqs T {\mathcal E}(T)&\leq&
C(\theta) (1+c_\alpha T) {\mathcal E}(0)+ C \int_{\Sigma_T}
(|\partial_tu|^2+|E_\tau|^2)\, dSdt+\theta \int_0^T {\mathcal
E}(t)\,dt
\\
&\leq& \big(C(\theta)  (1+c_\alpha T)+\theta T\big) {\mathcal E}(0)+
C \int_{\Sigma_T} (|\partial_tu|^2+|E_\tau|^2)\, dSdt, \forall
\theta>0,\eeqs reminding that the energy is non increasing. This is
the requested estimate (\ref{obsest}). \findem

\br {\rm Note that the last term of the estimate (\ref{sn3}) is zero
if $Q^*=Q^*_\nu$, but according to Remark \ref{rQnormal}, this
assumption is meaningless.} \er

\bt\label{tstabobs} Under the assumptions of the previous theorem
and if $c_\alpha$ is small enough,
 there exist
two positive constants $M$ and $\omega$ such that \be\label{so2}
 {\mathcal E}(t)\leq Me^{-\omega t} {\mathcal E}(0), \ee for all strong solution $(u,E,H)$ of
(\ref{max}). \et
\begin{remark}
{\rm The same method yields the same exponential stability result in
the case where $\eps, \mu$ are positive functions satisfying some
regularity and technical conditions.}
\end{remark}
\pf The estimate (\ref{obsest}) and Lemma
\ref{l2.1} yield
$$
T {\mathcal E}(T)\leq \big(C(\theta)  (1+ c_\alpha T)+\theta T\big)
{\mathcal E}(0)+ C ({\mathcal E}(0)-{\mathcal E}(T)),  \forall
\theta>0,$$ which may be equivalently written
$$
{\mathcal E}(T)\leq \frac{C(\theta)  (1+c_\alpha T)+\theta T}{C+T}
{\mathcal E}(0),\forall \theta>0.
$$
Now we choose  $\theta=\frac12$  and $c_\alpha\leq\frac{1}{4 \,
C(\frac12)}$, with this choice $\frac{C(\theta)  (1+ c_\alpha
T)+\theta T}{C+T}$ tends to $C(\frac12) c_\alpha+\frac12\leq
\frac34$ as $T$ goes to infinity. Therefore for $T$ large enough, we
have found $r\in (0,1)$ such that
$$
{\mathcal E}(T)\leq r {\mathcal E}(0).
$$
Since our system is invariant by translation, standard  arguments
about uniform stabilization of hyperbolic system (see for instance
\cite{RauchTaylor,ELN2}) yield the conclusion. \findem

The key point in the above proof is to estimate appropriately the
term $\int_{\Sigma_T}|u|^2\,dS dt$ in (\ref{star}). Indeed a rough
idea is to use the definition (\ref{2.1}) of the energy to get
\[
C\int_{\Sigma_T}|u|^2\,dS dt\leq \frac{2C}{A} \int_0^T{\mathcal
E}(t)\, dt\leq \frac{2C T}{A}  {\mathcal E}(0).\] Hence from the
previous proof we   obtain an exponential stability result only for
$A$ small enough (depending on a constant $C$ that is not known
explicitly, see nevertheless \cite{Alabau}). In order to prove the
stability result   for any positive $A$, we then need to estimate
$\int_{\Sigma_T}|u|^2\,dS dt$ in a different way. Its proof is based
on the use of a solution $z$ of a stationary problem (see
\cite{conrad:93,Alabau,HNicSene,nicaise} and below) such that $z=u$
on $\Gamma$. Multiplying the first identity of (\ref{max}) by $z$,
integrating by parts and using the second boundary condition in
(\ref{max}), the term $\int_{\Sigma_T}|u|^2\,dS dt$ naturally
appears. For standard problems (see
\cite{conrad:93,Alabau,HNicSene,nicaise}) this term is estimated
using elliptic regularity results on $z$. Here the specificity of
our piezoelectric system requires a more careful analysis. We start
with the stationary problem mentioned before.

\bl\label{l21.a} Let  $\left(u,E,H\right)$ be a  strong solution of
(\ref{max}). Then there exists $(z,\chi)\in H^1(\Om)^3\times
H^1_0(\Om)$ (depending on $t$) weak solution of \be\label{pbz}
\left\{
\begin{tabular}{ll}
& $\nabla (\sigma (z)-e^\top \nabla \chi)=0 $ in $\Om$,\vspace{2mm}
\\& $\ddiv (\eps \nabla \chi+e\gamma (z))=0 $ in $\Om$,\vspace{2mm}
\\
& $z=u, \chi=0$ on $\Gamma$.
\end{tabular}
\right. \ee Moreover there exists a positive constant $C$
(independent of $t$) such that \beq\label{l21estz} \int_\Om
|z|^2\,dx\leq C \int_\Gamma |u|^2\,dS\leq \frac{2C}{A}{\mathcal
E}(t),\\
\label{l21estz'} \int_\Om |\partial_t z|^2\,dx\leq C \int_\Gamma
|\partial_t u|^2\,dS\leq -C \partial_t {\mathcal E}(t). \eeq
 \el
\pf Inspired from \cite{conrad:93,Alabau,HNicSene,nicaise} for each
$t\geq 0$ we consider the weak solution $(z,\chi)$ (depending on
$t$) of (\ref{pbz}). This solution is characterized by   $z=w+u$
where $(w,\chi)\in \tilde{V}:=H^1_0(\Om)^3\times H^1_0(\Om)$ is the
unique solution of \be\label{pbwchi} \tilde{a}((w,\chi),
(w',\chi'))=-\tilde a((u,0), (w',\chi'), \forall (w',\chi')\in
\tilde{V}, \ee where
$$
\tilde{a}((w,\chi), (w',\chi'))= \int_\Om \{(\sigma (w)-e^\top
\nabla \chi): \gamma(w')+ (\eps \nabla \chi+e\gamma (w))\cdot \nabla
\chi'\} \,dx, \forall (w',\chi')\in V.$$ The above problem has a
unique solution since the bilinear form $\tilde a$ is coercive on
$V$ (consequence of Korn's inequality).

A direct consequence of (\ref{pbwchi}) is that
$$
\tilde{a}((z,\chi), (w',\chi'))= 0, \forall (w',\chi')\in \tilde{V}.
$$
By taking as test function $w'=w=z-u$ and $\chi'=\chi$, we find that
$$
\tilde{a}((z,\chi), (z,\chi))=\tilde{a}((z,\chi), (u,0)),
$$
which implies \be\label{l21pos} \int_\Om
\{\sigma(z):\gamma(u)-e^\top \nabla \chi: \gamma(u)\}\,dx=
\tilde{a}((z,\chi), (z,\chi))\geq 0. \ee

Note further that the coerciveness of $\tilde{a}$ leads to
$$
\|w\|_{1,\Om}+\|\chi\|_{1,\Om}\leq C \|u\|_{1,\Om},
$$
and then to \be\label{normeh1chiz}
 \|z\|_{1,\Om}+\|\chi\|_{1,\Om}\leq C \|u\|_{1,\Om}\leq C {\mathcal E}(t)^{1/2},\ee
 where $\|u\|_{s,\Om} = \|u\|_{H^s(\Om)}$.

Now we consider the adjoint problem: Find $(w^*,\chi^*)\in
\tilde{V}$ solution of \be\label{pbz*} \left\{
\begin{tabular}{ll}
& $\nabla (\sigma (w^*)+e^\top \nabla \chi^*)=z $ in
$\Om$,\vspace{2mm}
\\& $\ddiv (\eps \nabla \chi^*-e\gamma(w^*))=0 $ in $\Om$,\vspace{2mm}
\\
& $w^*=0, \chi^*=0$ on $\Gamma$,
\end{tabular}
\right. \ee which is the unique solution of  \be\label{pbwchi*}
\tilde{a}^*((w^*,\chi^*), (w',\chi'))=\int_\Om z \cdot w' \,dx,
\forall (w',\chi')\in V, \ee where
$$
\tilde{a}^*((w,\chi), (w',\chi'))= \int_\Om \{(\sigma (w)+e^\top
\nabla \chi): \gamma(w')+ (\eps \nabla \chi-e\gamma (w))\cdot \nabla
\chi'\} \,dx, \forall (w',\chi')\in V.$$ Again this problem has a
unique solution since the bilinear form $\tilde{a}^*$ is also
coercive on $\tilde{V}$. Since the system (\ref{pbz*}) is strongly
elliptic, we deduce that $(w^*,\chi^*)$ belongs to $H^2(\Om)^3\times
H^2(\Om)$ with the estimate (see Theorem 10.5 of
\cite{AgmonDouglisNirenberg64} or Theorem 4.5.3  of \cite{glc})
\be\label{pbwchi*est}\|w^*\|_{2,\Om}+\|\chi^*\|_{2,\Om} \leq C
\|z\|_{0,\Om}, \ee where here and below $C$ is a positive constant
that depends only on  $a_{ijkl},\eps,\mu, e_{ijk}$  and on $\Om$.

By using the differential equations from  (\ref{pbz*}), we may write
\beqs\int_\Om |z|^2 \,dx&=& \int_\Om \nabla (\sigma (w^*)+e^\top
\nabla \chi^*)\cdot z\,dx
\\
&=&\int_\Om \{\nabla  (\sigma (w^*)+e^\top \nabla \chi^*)\cdot z +
\ddiv (\eps \nabla \chi^*-e\gamma(w^*))\chi\}\, dx. \eeqs Applying
Green's formula we get \beqs\int_\Om |z|^2 \,dx&=& -\int_\Om
\{(\sigma (w^*)+e^\top \nabla \chi^*): \gamma(z) +   (\eps \nabla
\chi^*-e\gamma (w^*))\cdot \nabla \chi\}\, dx
\\
&+&\int_\Gamma (\sigma (w^*)+e^\top \nabla \chi^*) \nu \cdot z
\,dS \\
&=& -\int_\Om  \{(\sigma (z)-e^\top \nabla \chi): \gamma(w^*) +
(\eps \nabla \chi+e\gamma (z))\cdot \nabla \chi^*\}\, dx
\\
&+&\int_\Gamma (\sigma (w^*)+e^\top \nabla \chi^*) \nu \cdot z \,dS.
\eeqs Applying again Green's formula and reminding problem
(\ref{pbz}), we have found that
$$
\int_\Om |z|^2 \,dx=\int_\Gamma (\sigma (w^*)+e^\top \nabla \chi^*)
\nu \cdot u \,dS.
$$
By Cauchy-Schwarz's inequality and the estimate (\ref{pbwchi*est})
(with the help of a trace theorem), we obtain finally
\[ \int_\Om |z|^2\,dx\leq C \int_\Gamma |u|^2\,dS . \]
This proves (\ref{l21estz}) because $\frac{A}{2}\int_\Gamma
|u|^2\,dS \leq {\mathcal E}(t)$.

By deriving the system  (\ref{pbz}) in time, the estimate
(\ref{l21estz}) also shows that \[ \int_\Om |\partial_t z|^2\,dx\leq
C \int_\Gamma |\partial_t u|^2\,dS.\] This yields (\ref{l21estz'})
owing to the identity (\ref{2.2}). \findem

At this stage we need to exploit the fact that $\eps \nabla
\chi+e\gamma (z)$ is divergence free, hence it is the curl of
$\psi\in X_T(\Om)$, where \beqs X_T(\Om)=\{\phi\in H^1(\Om)^3: \ddiv
\psi=0 \hbox{ in }\Omega, \hbox{ and }  \psi\cdot \nu=0 \hbox{ on
}\Gamma\}. \eeqs More precisely we have the following result.

\bl\label{l21.b} Let  $\left(u,E,H\right)$ be a  strong solution of
(\ref{max}) and $(z,\chi)\in H^1(\Om)^3\times H^1_0(\Om)$  the weak
solution of (\ref{pbz}).  Then there exists $\psi\in X_T(\Om)$ such
that \be\label{egaliterotpsi} \eps \nabla \chi+e\gamma (z)= \rot
\psi, \ee with the estimates \beq\label{l21estpsi}
\|\psi\|_{0,\Om}^2\leq C \|u\|_{0,\Gamma}^2\leq  \frac{2C}{A}
{\mathcal E}(t), \\
\label{l21estpsiprime}
\|\partial_t\psi\|_{0,\Om}^2\leq C \|\partial_tu\|_{0,\Gamma}^2 \leq
-C
\partial_t {\mathcal E}(t). \eeq
where $C$ is a positive constant independent of $t$. \el \pf We
remark that   (see (\ref{pbz})) $\eps \nabla \chi+e\gamma (z)$ is
divergence free in $\Om$, hence as $\Om$ is simply connected, we
deduce (see Theorem I. 3.5 in \cite{GR})  that there exists $\psi\in
X_T(\Om)$ such that (\ref{egaliterotpsi}) holds  with the estimate
$$ \|\psi\|_{1,\Om}\leq C \, \|\eps \nabla \chi+e\gamma
(z)\|_{1,\Om}.
$$ Thanks to (\ref{normeh1chiz}), we get
 \be\label{normh1psi}
\|\psi\|_{1,\Om}  \leq C \|u\|_{1,\Om}\leq C {\mathcal
E}(t)^{1/2}.\ee

Let us finally consider the problem: find $\tilde\chi$ solution of
\be\label{pbchi} \left\{
\begin{tabular}{ll}
& $\rot \rot \tilde\chi =\psi $ in $\Om$,\vspace{2mm}
\\& $\ddiv  \tilde\chi=0 $ in $\Om$,\vspace{2mm}
\\
& $\tilde\chi\cdot \nu=0, \rot \tilde\chi\times \nu =0$ on $\Gamma$.
\end{tabular}
\right. \ee The variational formulation of this problem is: find
$\tilde\chi\in X_T(\Om)$ solution of \be\label{pbchivar}
b(\tilde\chi, \theta)=\int_\Om \psi \cdot \theta \,dx,   \forall
\theta\in H_T(\Om), \ee where
$$
b(\tilde\chi, \theta)= \int_\Om \{\rot \tilde\chi \rot \theta +\ddiv
\tilde\chi \ddiv \theta \}\,dx, \forall \tilde\chi, \theta\in
H_T(\Om),$$ and
$$
H_T(\Om)=\{\phi\in H^1(\Om)^3:  \phi\cdot \nu=0 \hbox{ on }\Gamma\}.
$$
It is well known (see for instance \cite{costabel:91}) that $b$ is
coercive on $H_T(\Om)$ and therefore problem (\ref{pbchivar}) is
well posed, its solution $\tilde \chi$ furthermore satisfies
(\ref{pbchi}) because $\psi$ is divergence free. Moreover as the
system $\rot \rot -\nabla\ddiv=-\Delta$ is strongly elliptic and the
boundary conditions in (\ref{pbchi}) cover this system, we get that
$\tilde\chi$ belongs to $H^2(\Om)^3$ with (see again Theorem 10.5 of
\cite{AgmonDouglisNirenberg64} or Theorem 4.5.3 of \cite{glc})
\be\label{estnormehetildechi}
 \|\tilde\chi\|_{2,\Om}
\leq C \|\psi\|_{0,\Om}. \ee Now as before we can write by using
Green's formula and the identity (\ref{egaliterotpsi}) \beqs
\|\psi\|_{0,\Om}^2&=&\int_\Om \psi\cdot\rot\rot \tilde\chi\,dx
\\
&=&\int_\Om \rot\psi\cdot\rot  \tilde\chi\,dx
\\
&=&\int_\Om (\eps \nabla \chi+e\gamma (z))\cdot\rot  \tilde\chi\,dx
\\
&=&-\int_\Om \nabla (e^\top \rot\tilde\chi)\cdot z\,dx +\int_\Gamma
(e^\top \rot\tilde \chi)\nu \cdot  z\,dS. \eeqs By the estimate
(\ref{estnormehetildechi}) and reminding that $z=u$ on $\Gamma$, we
obtain
$$
\|\psi\|_{0,\Om}\leq C (\|z\|_{0,\Om}+\|u\|_{0,\Gamma}).
$$
By the estimate (\ref{l21estz}), we arrive at \[
\|\psi\|_{0,\Om}^2\leq C \|u\|_{0,\Gamma}^2, \] and we conclude as
in the previous Lemma. \findem

\bl\label{l21} Let  $\left(u,E,H\right)$ be a  strong solution of
(\ref{max}). Then for all $\theta>0$ there exists a constant
$C(\theta)>0$ (which does not depend on $T$ but depends on $\theta$,
the domain and the coefficients $a_{ijkl},\eps,\mu, e_{ijk}, A$)
such that \be\label{estl21}
 \int_{\Sigma_T}|u|^2\,dS dt\leq C(\theta){\mathcal E}(0)
+\theta \int_0^T{\mathcal E}(t)\,dt. \ee \el \pf We multiply the
first identity of (\ref{max}) by $z\in H^1(\Om)^3$ from Lemma
\ref{l21.a} and integrate on $Q_T$ to get
$$
\int_{Q_T}z\cdot (\partial_t^2 u-\nabla \sigma (u,E))\,dxdt=0.
$$
By Green's formula we obtain
$$
\int_{Q_T}(z\cdot\partial_t^2 u+\sigma (u,E):\gamma(z))\,dxdt
-\int_{\Sigma_T} z\cdot(\sigma(u,E)\cdot \nu)\,dS dt=0.
$$
Using the second boundary condition in (\ref{max}) and the boundary
condition in (\ref{pbz}), we obtain \beqs A\int_{\Sigma_T} |u|^2\,dS
dt =-\int_{\Sigma_T} u \cdot (\partial_t u+ Q(E\times \nu))\,dS dt
-\int_{Q_T}(z\cdot\partial_t^2 u+\sigma (u,E):\gamma(z))\,dxdt.
\eeqs Owing to (\ref{l21pos}) we arrive at $$ A\int_{\Sigma_T}
|u|^2\,dS dt \leq
$$
\beqs -\int_{\Sigma_T} u \cdot (\partial_t u+ Q(E\times \nu))\,dS dt
-\int_{Q_T}(z\cdot \partial_t^2 u +e^\top
\nabla\chi:\gamma(u)-e^\top E:\gamma(z))\,dxdt. \eeqs By using the
identity $e\gamma(u)=D-\eps E$, we get \be\label{ideau2}
A\int_{\Sigma_T} |u|^2\,dS dt \leq-\int_{\Sigma_T} u \cdot
(\partial_t u+ Q(E\times \nu))\,dS dt  \ee
$$
-\int_{Q_T}\Big(z\cdot
\partial_t^2 u +\nabla \chi\cdot D-E\cdot
(e\gamma(z)+\eps\nabla\chi)\Big)\,dxdt. $$ We now transform the two
last terms of this identity, first by Green's formula in space, we
see that
$$
\int_{Q_T}\nabla \chi\cdot D \,dxdt= -\int_{Q_T} \chi\ddiv D \,dxdt
+\int_{\Sigma_T}  \chi  D\cdot\nu\,dS dt=0,
$$
since $D$ is divergence free and $\chi=0$ on $\Gamma$. On the other
hand, by the identity (\ref{egaliterotpsi}) we have
$$
\int_{Q_T}E\cdot (e\gamma(z)+\eps\nabla\chi)\,dxdt=
\int_{Q_T}E\cdot\rot \psi\,dxdt,
$$
and by Green's formula in space
$$
\int_{Q_T}E\cdot (e\gamma(z)+\eps\nabla\chi)\,dxdt= \int_{Q_T}\rot
E\cdot \psi\,dxdt+\int_{\Sigma_T}  (E\times \nu)\cdot \psi\,dS dt.
$$
Now reminding that $\mu \partial_t H=\rot E$ and using an
integration by parts in time, we arrive at
$$
\int_{Q_T}E\cdot (e\gamma(z)+\eps\nabla\chi)\,dxdt= \int_{Q_T} \mu
H\cdot \partial_t\psi\,dxdt+\int_\Om \mu H\cdot \psi
\,dx\Big|_0^T+\int_{\Sigma_T} (E\times \nu)\cdot \psi\,dS dt.
$$

In the same manner an integration by parts in time yields
$$
\int_{Q_T} z\cdot
\partial_t^2 u \,dxdt=-\int_{Q_T} \partial_tz\cdot
  \partial_tu\,dxdt+\int_\Om z\cdot
  \partial_tu\,dx\Big|_0^T
$$

These identities in (\ref{ideau2}) lead to

\beq \label{l21estu2} && A \, \int_{\Sigma_T} |u|^2\,dS dt
\leq-\int_{\Sigma_T} (u \cdot (\partial_tu+ Q(E\times \nu))+(E\times
\nu)\cdot \psi)\,dS dt
\\
 &&\hspace{2cm}+\int_{Q_T}(\partial_tz \partial_tu+ \mu H\cdot
\partial_t\psi)\,dxdt -\int_\Om z\cdot
  \partial_tu\,dx\Big|_0^T
  +\int_\Om \mu H\cdot \psi
\,dx\Big|_0^T. \nonumber
 \eeq

It remains to estimate each term of this right-hand side. For the
first term applying successively Cauchy-Schwarz's inequality,
Young's inequality and the identity (\ref{2.2}) we may write \beqs
\left|\int_{\Sigma_T} u \cdot (\partial_tu+ Q(E\times \nu))\,dS
dt\right| &\leq & \frac{A}{2} \int_{\Sigma_T} |u|^2\,dS dt
+\frac{C}{2A} \int_{\Sigma_T} \left(|\partial_tu|^2+|E\times \nu|^2
\right)\,dS dt
\\
&\leq & \frac{A}{2} \int_{\Sigma_T} |u|^2\,dS dt -\frac{C}{2A}
\int_0^T \partial_t {\mathcal E}(t) \,dt. \eeqs Since  the energy is
non-negative, we arrive at \be\label{l21term1} \left|\int_{\Sigma_T}
u \cdot
\partial_t u\,dS dt\right| \leq \frac{A}{2} \int_{\Sigma_T}
|u|^2\,dS dt +\frac{C}{2A}{\mathcal E}(0). \ee

For the second term by using Cauchy-Schwarz's inequality, Young's
inequality, a trace theorem, the estimate (\ref{normh1psi}) and
again the identity (\ref{2.2}) \beqs \left|\int_{\Sigma_T} (E\times
\nu)\cdot \psi\,dS dt\right| &\leq& \theta
\int_0^T\|\psi\|^2_{1,\Om}\,dt+\frac{C}{\theta} \int_{\Sigma_T}
|E\times \nu|^2\,dS dt
\\
&\leq& \theta \int_0^T{\mathcal E}(t)\,dt+\frac{C}{\theta}
\int_{\Sigma_T} |E\times \nu|^2\,dS dt
\\
&\leq& \theta \int_0^T{\mathcal E}(t)\,dt-\frac{C}{\theta} \int_0^T
\partial_t {\mathcal E}(t)
\, dt.
 \eeqs

As before   the energy being non-negative, we arrive at
\be\label{l21term1bis} \left|\int_{\Sigma_T} (E\times \nu)\cdot
\psi\,dS dt\right| \leq \theta \int_0^T{\mathcal
E}(t)\,dt+\frac{C}{\theta}{\mathcal E}(0). \ee

For the third term we use successively Cauchy-Schwarz's inequality,
Young's inequality, the estimate (\ref{l21estz'}) and the definition
of the energy to get for all $\theta>0$ \beqs \left|\int_{Q_T}
\partial_tz\cdot \partial_tu\,dxdt \right|&\leq&
\frac{1}{2\theta}\int_{Q_T}|\partial_tz|^2\,dxdt
+\frac{\theta}{2}\int_{Q_T}|\partial_tu|^2\,dxdt
\\
&\leq& -\frac{C}{2\theta}\int_0^T \partial_t {\mathcal E}(t) \,dt
+{\theta}\int_0^T {\mathcal E}(t) \,dt. \eeqs Again we get
\be\label{l21term2} \left|\int_{Q_T} \partial_tz\cdot
\partial_tu\,dxdt \right|\leq \frac{C}{\theta}{\mathcal
E}(0)+{\theta}\int_0^T {\mathcal E}(t) \,dt. \ee

As for the third term replacing the estimate (\ref{l21estz'}) by
(\ref{l21estpsiprime}) we get for the fourth term
\be\label{l21term2bis} \left|\int_{Q_T}\mu H\cdot
\partial_t\psi\,dxdt \right|\leq \frac{C}{\theta}{\mathcal
E}(0)+{\theta}\int_0^T {\mathcal E}(t) \,dt. \ee

For the fifth term the application of Cauchy-Schwarz's inequality,
the estimate (\ref{l21estz}) and the definition of the energy
directly gives
\be\label{l21term3.1}\left| \int_\Om z\cdot  \partial_tu
\,dx\Big|_0^T \right| \leq C({\mathcal E}(0)+{\mathcal E}(T))\leq 2C
{\mathcal E}(0) \ee since the energy is non-decreasing.

Similarly using (\ref{l21estpsi}) instead of (\ref{l21estz}), we
have

\be\label{l21term3}\left| \int_\Om \mu H\cdot \psi\,dx\Big|_0^T
\right| \leq C {\mathcal E}(0). \ee

The estimates (\ref{l21term1}) to (\ref{l21term3}) into the estimate
(\ref{l21estu2}) yield the conclusion. \findem

\end{document}